\documentclass[12pt,a4paper,reqno]{amsart}
\usepackage{amsfonts}
\usepackage{graphicx}   
\usepackage{graphics}
\usepackage{epsfig}
\usepackage{anysize}
\marginsize{2.5cm}{2.5cm}{2.5cm}{2.5cm}

\newcommand{\be}{\begin{equation}}
\newcommand{\bea}{\begin{eqnarray}}
\newcommand{\beas}{\begin{eqnarray*}}
\newcommand{\ee}{\end{equation}}
\newcommand{\eea}{\end{eqnarray}}
\newcommand{\eeas}{\end{eqnarray*}}
\def\begeq{\begin{equation}}
\def\endeq{\end{equation}}
\newcommand{\benu}{\begin{enumerate}}
\newcommand{\eenu}{\end{enumerate}}
\newcommand{\bit}{\begin{itemize}}
\newcommand{\eit}{\end{itemize}}

\def\Re{{\rm I\kern-.17em R}}

\def\({\left(}
\def\){\right)}

\begin{document}
\setcounter{page}{1}

\bigskip
\bigskip
\title[] {On a solution of one fuzzy logic problem}
\author[]{Shikhlinskaya R. Y. \dag}
\thanks{\dag AZ1148, Z.Khalilov, 23, Baku State University, Institute Applied Mathematics, Baku, Azerbaijan,\\ \noindent e-mail: reyhan$\_$sh@rambler.ru}

\begin{abstract}
In this paper defuzzification method of WABL is investigated, its
properties are analyzed. The WABL method is applied to some fuzzy
models. The package of applied programs is worked out on the base
of proposed algorithms. The obtained in the form of visual-
interactive graphs results are compared with knows ones.

\bigskip
 \noindent Keywords:  fuzzy number, defuzzification, averaging representative, waited average based, fuzzy logic.

\end{abstract}\maketitle

Fuzzy simulation nowadays is one of the intensive directions of
researches in control and decision making fields. Fuzzy
simulation, i.e. fuzzy theory-based program models enable to get
more adequate results in the field of control of engineering
systems. Their application field becomes more and more from year
to year: domestic electric equipments, cameras, video cameras,
industrial robots and others are converted to fuzzy logic based
working principles.

Sometimes to use the defuzzification instead of fuzzy parameters, i. e. exact representative is possible from the point of view of model's adequacy and this considerably simplifies the calculations. There are different defuzzification methods: COA(Center of Area), ÌÎÌ(Median of Maximum), OWA(Ordered Weighted Average),  WABL(Weighted Average Based on Levels) and ets. To use this or other defuzzification method depends on the statement of the problem.

 In the paper a new WABL defuzzification method-based logic deduction rules are given and applied to fuzzy regulators.

 We first consider WABL defuzzification method. Let the fuzzy number A be given by LR-description, i.e.

 \[A=\mathop{\bigcup }\limits_{\xi \in \left[0,1\right]}
\left(\xi ,A^{\xi } \right).\]

 Hear $A^{\xi } =\left[L_{A} \left(\xi \right),R_{A}
\left(\xi \right)\right]$ and $L_{A} \left(\xi \right)=\mu
_{\uparrow }^{-1} \left(\xi \right),\, \, \, \, \, R_{A} \left(\xi
\right)=\mu _{\downarrow }^{-1} \left(\xi \right).$ $\mu
_{\uparrow }^{-1} $ and $\mu _{\downarrow }^{-1} $ define left and
right sides of a fuzzy number correspondingly.

\textbf{Definition.}When we say averaged representative  of fuzzy
number or WABL aggregation we mean the following quantity:

\begin{equation} \label{GrindEQ__1_}
I_{W} (A)=c_{L} \int _{0}^{1}L_{A} (\xi )p(\xi )d\xi +c_{R} \int
_{0}^{1}R_{A} (\xi )p(\xi )d\xi
\end{equation}

Where $c_{L} $ and $c_{R} $ indicate the importance degree of the left and right hand sides, respectively, by operating on fuzzy members and

\[\begin{array}{l} {{\rm c}_{L} \ge 0,\, \, \, \, \, \, c_{R} \ge 0:}  {c_{L} +c_{R} =1} \end{array}\]

\noindent $p\left(\xi \right)$ is the distribution function of importance degrees.

\noindent

\noindent $P:\left[0,1\right]\to E_{+} \equiv \left[0,+\infty
\right]$  and  $\int _{0}^{1}p\left(\xi \right)d\xi =1 $

We can show that a general form of the function $p\left(\xi \right)$ is as follows:

\noindent

\begin{equation} \label{GrindEQ__2_}
p(\xi )=k\xi ^{k-1}
\end{equation}
here $k>0$.

\underbar{The properties of WABL defuzzification method}.

1. WABL method is more general. For some class of fuzzy numbers we can select the parameters $c_{L} $ , $c_{R} $ , and  $p\left(\xi \right)$  so that the result obtained by WABL defuzzification method coincides with the result obtained by an other defuzzification methods.

2. WABL defuzzification method is additive,   i.e.

\[I_{W} (A+B)=I_{W} (A)+I_{W} (B)\]
This method enables to take into account the expert's strategy in a model.

 3. If the result of defuzzification doesn't satisfy the experts by changing $c_{L} $ , $c_{R} $ , and  $p\left(\xi \right)$  parameters he may change defuzzification number or slide it arbitrarily to the left or right hand side.

 Now give the stages of logic deduction rules based on WABL defuzzification method.

1.  Organization of rules base for logic deduction rules;

2. Fuzzification of input and output variables, i.e. construction
of fuzzy
 sets $\mu _{k_{i} } (x)$, $i=\overline{1,n}$  and $\mu
_{l_{i} } (x)$, $j=\overline{1,m}$  describing the $k_{i} $  input
fuzzy and $l_{i} $ output fuzzy states;

 3. Calculation of validity degree of entrance of $x=\overline{x}$  exact input estimate to the fuzzy input estimates set $\mu _{k_{i} } (x)$, $i=\overline{1,n}$. Denote this validity degrees by  $\overline{k_{i} }$: $\overline{k_{i} }=\mu _{k_{i} } (\overline{x})$, $i=\overline{1,n}$;

 4. Defuzzification of $l_{i} $ output fuzzy subsets by WABL method (determination of  $I_{W} (l_{j} )$, $j=\overline{1,m}$;

 5. Addition of the exact input value of defuzzification results with the coefficients equal to validity degree;

\[I_{W} =\sum _{i=1}^{n}\sum _{j=1}^{m}\overline{k_{i} }I_{W} (l_{j} )  \]

Example. Let's see the functioning fuzzy deduction rules in expert system controlling a fan of room conditioner (see 1, p. 262) . Conditioner's work keeps optimal temperature in the room .

Assume that by changing rotation speed of  the fan . we can change the temperature of the room. let's construct  the working algorithm of the conditioner step by step. We first give rules base for the working algorithm. Notice that the organization of  rules base in a subjective process and depends on the desire of the expert .

Rules base :

\noindent   1) If the temperature of the room is \underbar{lover}, the rotation speed of a fan is \underbar{lover;}

\noindent   2) If the temperature of the room is \underbar{middle}, the rotation speed of a fan is \underbar{middle;}

\noindent   3) If the temperature of the room is \underbar{higher}, the rotation speed of a fan is \underbar{higher;}

As  we see there are six fuzzy terms sets in rules base. Three of
them:   "lower temperature", ''middle temperature `` and  ''
higher temperature `` are  input fuzzy parameters; the other ones
:''lover speed `` , ''middle speed ``,''higher speed `` are output
fuzzy parameters.

 We are to define membership functions for fuzzy subsets determined in fuzzy input value set of temperature and in fuzzy output values set of rotation speed of the fan in defuzzification step of fuzzy variables.

 Let's define membership function for the  ''lover``, ''middle`` and ''higher``  fuzzy subsets of temperature.

\begin{equation} \label{GrindEQ__3_}
\mu _{t_{lower} } (t)=\left\{\begin{array}{l} {1,\, \, \, \, \, \, \, \, \, \, \, \, \, \, \, \, \, \, \, \, \, \, \, t\le 12} \\ {\frac{20-t}{8} ,\, \, \, \, \, \, \, \, \, 12<t<20} \\ {0,\, \, \, \, \, \, \, \, \, \, \, \, \, \, \, \, \, \, \, \, \, t\ge 20} \end{array}\right.
\end{equation}

\begin{equation} \label{GrindEQ__4_}
\mu _{t_{middle} } (t)=\left\{\begin{array}{l} {0,\, \, \, \, \, \, \, \, \, \, \, \, \, \, \, \, \, \, \, \, \, \, t<12\, \, ve\, \, \, t>30} \\ {\frac{t-12}{8} ,\, \, \, \, \, \, \, \, \, \, 12\le t<20} \\ {\frac{30-t}{10} ,\, \, \, \, \, \, \, \, \, \, 20\le t\le 30} \end{array}\right.
\end{equation}

\begin{equation} \label{GrindEQ__5_}
\mu _{t_{higher} } (t)=\left\{\begin{array}{l} {0,\, \, \, \, \, \, \, \, \, \, \, \, \, \, \, \, \, \, \, \, \, \, 0<t<20} \\ {\frac{t-20}{10} ,\, \, \, \, \, \, \, \, \, \, 20\le t<30} \\ {1,\, \, \, \, \, \, \, \, \, \, \, \, \, \, \, \, \, \, \, \, \, \, 30\le t<60} \end{array}\right.
\end{equation}

 Now let's define fuzzy subsets for output variables. Assume that the rotation speed of the fan varies from
 0 and 1000 rot/min. then we'll define for this speed the ''lower``, ''middle`` and ''higher``
  fuzzy subsets and their membership functions as follows:

\begin{equation} \label{GrindEQ__6_}
\mu _{v_{lower} } (t)=\left\{\begin{array}{l} {1,\, \, \, \, \,
\, \, \, \, \, \, \, \, \, \, \, \, \, \, \, \, \, 0<v<200} \\
{\frac{400-v}{200} ,\, \, \, \, \, \, \, 200\le v\le 400} \\ {0,\,
\, \, \, \, \, \, \, \, \, \, \, \, \, \, \, \, \, \, \, \, \,
400<v<1000} \end{array}\right.
\end{equation}

\begin{equation} \label{GrindEQ__7_}
\mu _{v_middle}(t)=\left\{\begin{array}{l} {0,\, \, \, \, \, \, \,
\, \, \, \, \, \, \, \, \, \, \, \, \, \, \, 600<v<200}
\\ {\frac{v-200}{200} ,\, \, \, \, \, \, \, \, \, \, 200\le v<400}
\\ {\frac{600-v}{200} ,\, \, \, \, \, \, \, \, \, \, 400\le v\le
600} \end{array}\right.
\end{equation}

\begin{equation} \label{GrindEQ__8_}
\mu _{v_{higher} } (t)=\left\{\begin{array}{l} {0,\, \, \, \, \, \, \, \, \, \, \, \, \, \, \, \, \, \, \, \, \, \, 0<v<400} \\ {\frac{v-400}{200} ,\, \, \, \, \, \, \, \, 400\le v<600} \\ {1,\, \, \, \, \, \, \, \, \, \, \, \, \, \, \, \, \, \, \, \, \, \, \, 600\le v<1000} \end{array}\right.
\end{equation}

\noindent       3. Let's see the process of definition of rotation speed of a fan depending a temperature of the room by fuzzy expert system . Assume the room's temperature is $t=22^{0} C$  .

\noindent       The first calculate the degree of  $t=22^{0} C$   value  to  fuzzy subset \eqref{GrindEQ__1_} -- \eqref{GrindEQ__3_}:

\noindent

\[\mu _{t_{lower} } (22)=0\]

\begin{equation} \label{GrindEQ__9_}
\mu _{t_{middle} } (22)=0,8
\end{equation}

\[\mu _{t_{higher} } (22)=0,2\]

\noindent In this step      defuzzificate fuzzy output subset \eqref{GrindEQ__6_}  - \eqref{GrindEQ__8_} by WABL method.

Let's write the increasing and diminishing parts in fuzzy set \eqref{GrindEQ__7_} according to levels:

\[L_{\nu _{lower} } (\xi )=0,    R_{v_{lower} } (\xi )=-200\xi +400\]

 Take them into account in \eqref{GrindEQ__1_}:

\noindent

\[I_{W} (v_{lower} )=(1-c_{L} )\int _{0}^{1}(-200\xi +400)(k+1)\xi ^{k} d\xi \, = \, -200c_{L} \frac{k+3}{k+2} -\]

\begin{equation} \label{GrindEQ__10_}
\, \, \, \, \, \, \, \, \, \, \, \, \, \, \, \, \, \, \, \, \, \, \, \, \, \, \, \, \, \, \, \, \, \, \, \, \, \, \, \, \, \, \, \, \, -\, 200\frac{k+1}{k+2} +400.
\end{equation}

For fuzzy set \eqref{GrindEQ__7_}-\eqref{GrindEQ__8_} :

\[L_{v_{middle} } (\xi )=200\xi +200,   R_{v_{_{middle} } } (\xi )=-200\xi +600,\]

\begin{equation} \label{GrindEQ__11_}
\begin{array}{l} {I_{W} (v_{middle} )=c_{L} \int _{0}^{1}(200\xi +200)(k+1)
\xi ^{k} d\xi +(1-c_{L} )\int _{0}^{1}(-200\xi +  } \\
\\ { +600)(k+1)\xi ^{k} d\xi =-400c_{L}
\frac{1}{k+2} -200\frac{k+1}{k+2} +600,} \end{array}
\end{equation}

\[L_{v_{B} higher} (\xi )=200\xi +400,   R_{v_{higher} } (\xi )=1000,\]

\begin{equation} \label{GrindEQ__12_}
\begin{array}{l} {I_{W} (v_{higher} )=c_{L} \int _{0}^{1}(200\xi +400)(k+1)\xi d\xi +(1-c_{L} )
\int _{0}^{1}1000(k+1)\xi ^{k} d\xi =  } \\
\\{ =-200c_{L} \frac{2k+5}{k+2} +1000} \end{array}
\end{equation}

        Now multiply the averaged representative of input fuzzy subsets  $l_{i} $ by corresponding constants equal to  $k_{i} $ and summarize them:

\begin{equation} \label{GrindEQ__13_}
I_{W} =0\cdot I_{W} (l_{lower} )+0.8\cdot I_{W} (l_{middle} )+0.2\cdot I_{W} (l_{higher} )=0.8\cdot 400+0.2\cdot 766=473.2
\end{equation}

So if temperature will be $ 22^{0} C$ the expert will adopt the
rotation speed of

\noindent a fan  473.2 rot/min.

Each time when we result of logical deduction is required the above

\noindent mentioned operations are carried out by the same order. In dynamic control systems this sequence is performed periodically.

    It each time we'll mark the pair $(t,v)$  in coordinate surface we'll obtain a graph of the dependence of rotation speed.

It was proved that fuzzy logic conditioners provide the less vibration of temperature and saves much electric power.
The work of fuzzy logic based systems in dynamic processes is more effective than work of
  control circuits working on an ordinary thermostat.

 The other expert system may determine the rotation speed by another way.
 But by changing the parameters $c_{L_{} } ,c_{R} $ and $p(\xi )$  in WABL defuzzification  formula. We can arbitrarily change the rotation speed of a fan  within  its definition domain. The comparing the results of the work of our expert system and the expert system given in [1] with different values of input parameters is shown on the p. 3-4. The optimal temperature to what the temperature of the room tents as a result  of air conditioner working is defined by itself system. In our case changing the values of $c_{L_{} } $and $c_{R} $. The optimal temperature of the room may be changed.

\bigskip\bigskip\bigskip

                                                 \begin{center}
                                                 \textbf{Referenses}\end{center}

\noindent

1. Data based. The intellectual processing of the information. /
V.V. Korneev, M., Nolij, 2001, p.255-290.

2.     Leonenkov À. The fuzzy modeling in MATLAB and fuzzy TECH.
SP v.: BÕV, Peterburg, 2003, p.73.

3.   Nasibov E.N. The methods for information processing in the
designed making problems. Baku: «Elm», 2000, p. 260.

4. Nasibov E.N., Shikhlinskaya R.Y.  Set-up of WABL-aggregation
parameters for finding the gravity centre of fuzzy number of
triangle form // Translations of Azerbaijan National Academy of
Sciences  (2003) No.2.

5.  Nasibov E.N., Shikhlinskaya R.Y. Set-up of WABL-aggregation
parameters for finding the gravity centre of fuzzy number of
polynomial form // Automatic Control and Computer Sciences 37
(2003) No.6. (in English).

\end{document}